\documentclass[secthm]{elsart}
\usepackage{times}
\usepackage{amsfonts}
\usepackage{amsmath}
\usepackage{amscd}
\usepackage{amssymb}
\usepackage{latexsym}
\usepackage{natbib}

\newcommand\preprintnote {preprint on \myhomepage}
\newcommand\myhomepage{http://www.math.ohio-state.edu/\-\~{}schoutens}

\newcommand\OSU{\address{Department of Mathematics\\
100 Math Tower\\
Ohio State University\\
Columbus, OH 43210 (USA)}
}

\newcommand{\emptyprop}{q}

\newcommand \ann[2]{\operatorname{Ann}_{#1}(#2)}
 
\newcommand \exactseq [5]{0\to{#1}\:\map{#2}\:{#3}\:\map{#4}\:{#5}\to0}
\newcommand \Exactseq [3]{0\to {#1}\to {#2}\to {#3}\to 0}

\newcommand \id{\mathfrak a}

\newcommand \iso{\cong}

\newcommand \map[1]{{\newcommand{\tmpprop}{#1q}  \if\tmpprop\emptyprop \to\else \xrightarrow{{\phantom{i}{#1}\phantom{i}}}\fi}} 
\newcommand \maxim{\mathfrak m}

\newcommand \pr{\mathfrak p}

\newcommand \rij[2]{(#1_1,\dots,#1_{#2})}

\newcommand \tensor{\otimes}
\newcommand \tor[4]{\operatorname{Tor}^{#1}_{#2}(#3,#4)}
\newcommand \op\operatorname

\newcommand{\name}[1]{#1}

\newcommand \ch{characteristic}
\newcommand \homo{homomorphism}
\newcommand \CM{Coh\-en-Mac\-au\-lay}
\renewcommand\iff{if, and only if,}

\newcommand  \BCM[1]{\mathcal B(#1)}

\begin{document}

\begin{frontmatter}
\title{On the vanishing of Tor of the absolute integral closure}
\author{Hans Schoutens}
\thanks{Partially supported by a  grant from the National Science Foundation.}
\date{02.04.2003}
\OSU
\ead{schoutens@math.ohio-state.edu}
\ead[url]{www.math.ohio-state.edu/~schoutens}

\begin{abstract}   
Let $R$ be an excellent local domain of positive \ch\ with residue field $k$ and let $R^+$ be its absolute integral closure. If $\tor R1{R^+}k$ vanishes, then $R$ is \CM, normal, F-rational and F-pure. If $R$ has at most an isolated singularity or has dimension at most two, then $R$ is regular.
\end{abstract}



\end{frontmatter}

\section{Introduction}

Recall that the \emph{absolute integral closure} $A^+$ is defined for an arbitrary domain $A$ as the integral closure of $A$ inside an algebraic closure of the field of fractions of $A$. A key property of the absolute integral closure was discovered by \name{Hochster} and \name{Huneke} in \citep{HHbigCM}: for $R$ an excellent local domain of positive \ch, $R^+$ is a \emph{balanced big \CM\ algebra}, that is to say,  any system of parameters on $R$ is an $R^+$-regular sequence. It is well-known that this implies that an excellent local domain $R$ of positive \ch\  is regular \iff\  $R\to R^+$ is  flat. Indeed, the direct implication follows since $R^+$ is a balanced big \CM\ algebra of finite projective dimension (use for instance \citep[Theorem IV.1]{SchFPD}) and the converse follows since $R\to R^+$ and $R^{1/p}\to R^+$ are isomorphic whence both faithfully flat, implying that $R\to R^{1/p}$ is flat, and therefore, by \name{Kunz}'s Theorem, that $R$ is regular (here $R^{1/p}$ denotes the extension of $R$ obtained by adding all $p$-th roots of element of $R$; for more details see \citep[Theorem 9.1 and Exercise 8.8]{HuTC}).

In \citep[Exercise 8.8]{HuTC}, \name{Huneke} points out that it is not known whether the weaker condition that  all \emph{Betti numbers} of $R^+$ vanish, that is to say, that all $\tor Rn{R^+}k$ vanish for $n\geq 1$, already implies that $R$ is regular. It is not hard to see, using that $R^+$ is a big \CM\ algebra, that this is equivalent with requiring that only $\tor R1{R^+}k$ vanishes. The main result of this paper is then the following positive solution for isolated singularities.

\begin{thm}\label{T:main}
Let $(R,\maxim)$ be an excellent  local domain   of positive \ch\ with residue field $k$.  Suppose $R$ has either an isolated singularity or  has dimension at most two. If $\tor R1{R^+}k=0$, then $R$ is regular.
\end{thm}

For arbitrary domains,  we obtain at least the following.

\begin{thm}\label{T:pure}
Let $(R,\maxim)$ be an excellent  local domain   of positive \ch\ and suppose $\tor R1{R^+}k=0$, where $k$ is the  residue field of $R$.  Then $R$ is F-rational and F-pure and any finite extension of $R$ is split. In particular, $R$ is normal, \CM\ and pseudo-rational, and every ideal is Frobenius-closed and plus-closed.
\end{thm}

We have some more precise information on the vanishing of certain \emph{Tor}'s in terms of the  singular locus  of $R$.

\begin{thm}\label{T:sing}
Let $(R,\maxim)$ be an excellent  local domain   of positive \ch\ and let $\id$ be an ideal defining the singular locus of $R$ (e.g., $\id$ is the \emph{Jacobian} of $R$).  If $\tor R1{R^+}k=0$, where $k$ is the  residue field of $R$, then $\tor Rn{R^+}M=0$ for all $n\geq 1$ and all finitely generated $R$-modules $M$ for which $M/\id M$ has finite length.
\end{thm}

The key observation in obtaining all these results, is that, in general, the vanishing of $\tor R1Sk$ implies that $R\to S$ is \emph{cyclically pure} (or \emph{ideal-pure}), meaning that $IS\cap R=I$, for all ideals $I$ of $R$. This is explained in Section~\ref{s:tor}. Moreover, \name{Hochster} has shown in \citep{HoPure} that for an excellent normal domain, cyclic purity is equivalent with purity, so that Theorem~\ref{T:pure} follows from some results of \name{Smith} on F-rational rings and plus-closure of parameter ideals. To prove Theorem~\ref{T:main}, we need a result  from \citep{SchFPD}: if  the first Betti number of a module over an isolated singularity vanishes, then the module has finite projective dimension.  Now, the argument which proofs that $R\to R^+$ is flat when $R$ is regular, yields the same conclusion under the weaker assumption that $R^+$ has finite projective dimension. This proves also the two-dimensional case, since we know already that $R$ is normal. 

Balanced big \CM\ algebras in \ch\ zero exist by the work of \name{Hochster-Huneke}, basically by a lifting procedure due to \name{Hochster}. However, the balanced big \CM\ algebras obtained in \citep{HHbigCM} are not canonically defined. In \citep{SchBCM}, I give an alternative but canonical construction $\BCM R$ of a balanced big \CM\ algebra  for a $\mathbb C$-affine local domain $R$ using ultraproducts and the absolute integral closure in positive \ch. It will follow form the present results that if $\tor R1{\BCM R}k=0$, where $k$ is the residue field of $R$, then $R$ is regular provided $R$ has an isolated singularity or has dimension at most two (moreover,without these additional assumptions, $R$ has at most rational singularities). This is the more interesting because it is not clear whether in general flatness of $R\to \BCM R$ implies regularity of $R$.

\section{Vanishing of  Betti numbers and cyclic purity}\label{s:tor}

We derive a simple criterion for a local ring \homo\ to be cyclically pure. We start with an easy lemma, the proof of which is included for sake of completeness. 

\begin{lem}\label{L:ss}
Let $A$ be a ring, $\id$ an ideal in $A$ and $ M$ and $N$ two $A$-modules. If $\id N=0$ and $\tor A1 MN=0$, then $\tor {A/\id}1{ M/\id M} N=0$.
\end{lem}
\begin{pf}
One can  derive this by aid of spectral sequences, but the following argument is more direct. Put $\bar A:=A/\id$. Since $ N$ is an $\bar A$-module, we can choose an exact sequence of $\bar A$-modules
	\begin{equation*}
	\Exactseq{\bar H}{\bar F} N
	\end{equation*}
with $\bar F$ a free $\bar A$-module. Tensoring with the $\bar A$-module $\bar M:= M/\id M$, we get an exact sequence
	\begin{equation*}
	0 \to \tor {\bar A}1{\bar M}{ N}\to \bar M\tensor_{\bar A}\bar H \to \bar M\tensor_{\bar A}\bar F.
	\end{equation*}
Since the last two modules are equal to $M\tensor_A \bar H$ and $M\tensor_A \bar F$ respectively and since $\tor A1 M N=0$, the last morphism in this exact sequence is injective. Therefore, $\tor {\bar A}1{\bar M}{ N}=0$, as required.
\qed\end{pf}

\begin{thm}\label{T:tor}
Let $(R,\maxim)$ be a Noetherian local ring with residue field $k$ and let $S$ be an arbitrary $R$-algebra. If $\tor R1S k=0$ and $\maxim S\neq S$, then $R\to S$ is cyclically pure.
\end{thm}
\begin{pf}
Since $\tor R1Sk$ vanishes,  so does $\tor{R/\mathfrak n}1{S/\mathfrak nS}k$ by Lemma~\ref{L:ss}, for every $\maxim$-primary ideal $\mathfrak n$. By the Local Flatness Criterion (see \citep[Theorem 22.3]{Mats}) applied to the Artinian local ring $R/\mathfrak n$,  the base change $R/\mathfrak n\to S/\mathfrak nS$ is flat, whence faithfully flat, since $\maxim S\neq S$. In particular, this base change is injective, showing that $\mathfrak nS\cap R=\mathfrak n$. Since every ideal is the intersection of $\maxim$-primary ideals by Krull's Intersection Theorem, the assertion follows.
\qed\end{pf}

\begin{rem}
Note that with notation from the Theorem, we have that the induced map of affine schemes $\op{Spec}S\to \op{Spec}R$ is surjective, since the \emph{fiber rings} $S_\pr/\pr S_\pr$ are non-zero.
\end{rem}

The following lemma shows that for a local \CM\ ring, the  vanishing of some Betti number of a  big \CM\ algebra is equivalent with the vanishing of all of its Betti numbers. 

\begin{lem}\label{L:rig}
If $(R,\maxim)$ is a local \CM\ ring with residue field $k$ and if $S$ is a big \CM\ $R$-algebra, such that $\tor RjSk=0$ for some $j\geq 1$, then $\tor RnSk=0$, for all $n\geq 1$.
\end{lem}
\begin{pf}
Let $\mathbf x$ be a maximal $R$-regular sequence which is also $S$-regular. Put $I:=\mathbf xR$. Since  $\tor RjSk$ vanishes, so does $\tor {R/I}j{S/IS}k$, so that $S/IS$ has finite flat dimension over $R/I$ by the Local Flatness Criterion. However, since the finitistic weak dimension is at most the dimension of a ring by \citep[Theorem 2.4]{ABII}, it follows that $S/IS$ is flat over $R/I$. Therefore, $0=\tor {R/I}n{S/IS}k=\tor RnSk$, for all $n\geq 1$.
\qed\end{pf}

Therefore, below, we may replace everywhere the condition that $\tor R1Sk=0$ by the weaker condition that some $\tor RjSk=0$, provided we also assume that $R$ is \CM. In fact, if $j$ is either $1$ or $2$, we do not need to assume that $R$ is \CM, since this then holds automatically.

\begin{prop}
If $(R,\maxim)$ is a Noetherian local   ring with residue field $k$ and if $S$ is a big \CM\ $R$-algebra, such that either $\tor R1Sk$ or $\tor R2Sk$ vanishes, then $R$ is \CM.
\end{prop}
\begin{pf}
I claim that $IS\cap R=I$, for some parameter ideal of $R$. By a standard argument, it then follows  that $R$ is \CM\ (see for instance the argument in \citep[Theorem 4.2]{SchBCM}). For $j=1$, we can use Lemma~\ref{L:ss} to conclude that $\tor {R/I}1{S/IS}k=0$, so that by the argument above, $R/I\to S/IS$ is faithfully flat. For $j=2$, we reason as follows. Let
	\begin{equation*}
	\Exactseq MFS
	\end{equation*}
be a short exact sequence with $F$ free. It follows that $\tor R1Mk$ is equal to $\tor R2Sk$, whence is zero. Therefore, by the same argument as before, $M/IM$ is flat over $R/I$. On the other hand, since we may choose $I$ so that it is generated by an $S$-regular sequence, we get that $\tor R1S{R/I}=0$ (indeed, the canonical morphism $I\tensor M\to IM$ is easily seen to be injective). Hence we get an exact sequence
	\begin{equation*}
	\Exactseq {M/IM}{F/IF}{S/IS}
	\end{equation*}
showing that $S/IS$ has finite flat dimension, whence is flat, since $R/I$ is Artinian.
\qed\end{pf} 

Is there a counterexample in which some $\tor RjSk$ vanishes for some big \CM\ algebra $S$ and some $j>2$, without $R$ being \CM?

\section{Proofs}

Recall that an excellent local ring of positive \ch\ is called \emph{F-rational}, if some ideal generated by a system of parameters  is tightly closed. It is well-known that an F-rational ring is \CM\ and normal, whence in particular a domain (\citep[Theorem 4.2]{HuTC}). By \citep[Theorem 3.1]{SmFrat}, an F-rational ring is pseudo-rational. We say that $R$ is \emph{F-pure}, if $R\to R^{1/p}$ is pure (or, equivalently, if the Frobenius is pure). 

\subsubsection*{Proof of Theorem~\ref{T:pure}}
Suppose $R$ is as in the statement of the theorem, so that in particuar $\tor R1{R^+}k$ vanishes. By Theorem~\ref{T:tor}, the embedding $R\to R^+$ is cyclically pure. Let $I$ be an ideal of $R$ generated by a system of parameters. It follows that $IR^+\cap R=I$. By the result in \citep{SmParId} that plus-closure and tight closure agree on parameter ideals, we get that $I$ is tightly closed. Hence $R$ is F-rational and therefore normal, \CM\ and pseudo-rational, by our previous observations (in fact, cyclic purity of $R\to R^+$ together with the fact that $R^+$ is a balanced big \CM\ algebra shows immediately that $R$ is \CM\ and normal). Since $R$ is normal, it follows  from \citep{HoPure} that $R\to R^+$ is pure. Let $R\subset  S$ be a finite extension. In order to show that this is split, we may factor out a minimal prime of $S$ and hence assume that $S$ is a domain. So $R\subset S$ extends to the pure map $R\to R^+$ and hence is itself pure. Since a pure map with finitely generated cokernel is split (\citep[Theorem 7.14]{Mats}), we showed that any finite extension splits. Finally, since by our previous argument $R\to R^{1/p}$ splits, $R$ is F-pure.
\qed

\subsubsection*{Proof of  Theorem~\ref{T:main}}
The vanishing of $\tor R1{R^+}k$ implies that $R$ is \CM\ by Theorem~\ref{T:pure}. Since $R^+$ is a balanced big \CM\ algebra and since $R$ has an isolated singularity, we get from \citep[Theorem IV.1]{SchFPD} that $R\to R^+$ is flat. As already observed, this implies that $R$ is regular. If $R$ has dimension at most $2$, then by Theorem~\ref{T:pure}, it is normal and therefore has an isolated singularity, so that the previous argument applies.
\qed

Recall that by the argument at the end of the previous section, the vanishing of some $\tor Rj{R^+}k$ implies already that $R$ is regular, if apart from being an isolated singularity, we also assume that $R$ is \CM, when $j\geq3$. In order to derive a regularity criterion from  Theorem~\ref{T:main}, we need a lemma on flatness over Artinian local Gorenstein rings of embedding dimension one.

\begin{lem}\label{L:flat1}
Let $(A,\maxim)$ be an Artinian local ring of embedding dimension one and let $M$ be an arbitrary $A$-module. Then $M$ is $A$-flat \iff\ $\ann MI=\maxim M$, where $I$ denotes the socle of $A$, that is to say, $I=\ann A\maxim$. 
\end{lem}
\begin{pf}
By assumption $\maxim=xA$, for some $x\in A$. It follows that the socle $I$ of $A$ is equal to $x^{e-1}A$, where $e$ is the smallest integer for which $x^e=0$. I claim that $\ann M{x^{e-i}}=x^iM$, for all $i$. We will induct on $i$, where the case $i=1$ is just our assumption. For $i>1$, let $\mu\in M$ be such that $x^{e-i}\mu=0$. Therefore, $x^{e-i+1}\mu=0$, so that by our induction hypothesis, $\mu\in x^{i-1}M$, say, $\mu=x^{i-1}\nu$. Since $0=x^{e-i}\mu=x^{e-1}\nu$, we get $\nu\in xM$ whence $\mu\in x^iM$, as required.

Flatness now follows by the Local Flatness Criterion \citep[Theorem 22.3]{Mats}. Indeed, it suffices to show that $A/xA\to M/xM$ is flat and $xA\tensor M\iso xM$. The first assertion is immediate since $A/xA$ is a field. For the second assertion, observe that $xA\iso A/x^{e-1}A$ and by what we just proved $xM\iso M/\ann Mx\iso M/x^{e-1}M$. It follows that $xA\tensor M$ is isomorphic with $xM$, as required.
\qed\end{pf}

\begin{cor}
Let $(R,\maxim)$ be a $d$-dimensional excellent  local \CM\ domain   of positive \ch.  Suppose that there exists an  ideal $I$ in $R$ generated by a regular sequence such that $\maxim/I$ is a cyclic module. Suppose also that $R$ has either an isolated singularity or that $d\leq 2$. If for each finite extension domain $R\subset S$, we can find a finite extension $S\subset T$, such that
	\begin{equation}\label{eq:socin}
	(IS:_S (I:_R\maxim)S)\subset \maxim T,
	\end{equation}
then $R$ is regular.
\end{cor}
\begin{pf}
Let $\rij xi$ be the regular sequence generating $I$ and write $\maxim=I+xR$. If $i<d$ then necessary $i=d-1$ and $\maxim$ is generated by $d$ elements, so $R$ is regular. Hence assume $i=d$, that is to say, $I$ is $\maxim$-primary. It follows that $\overline R:=R/I$ is an Artinian local ring with maximal ideal $x\overline R$.  Let $e$ be the smallest integer for which $x^e\in I$. Hence the socle of $\overline R$ is $x^{e-1}\overline R$. Let $\overline{R^+}:=R^+/IR^+$.  I claim that 
	\begin{equation*}
	\ann {\overline{R^+}}{x^{e-1}} = x{\overline{R^+}}.
	\end{equation*}
Assuming the claim, Lemma~\ref{L:flat1} yields that $\overline{R^+}$ is $\overline R$-flat. Therefore, if $k$ is the residue field of $R$, then $\tor {\overline R}1{\overline{R^+}}k=0$. But $\rij xd$ is both $R$-regular and $R^+$-regular, so that $\tor R1{R^+}k=0$. Regularity of $R$ then follows from Theorem~\ref{T:main}.

To prove the claim, one inclusion is clear, so assume that $a\in R^+$ is such that $ax^{e-1}\in IR^+$. Choose a finite extension $R\subset S\subset R^+$ containing $a$ and such that we already have a relation $ax^{e-1}\in IS$. By assumption, we can find a finite extension $T$ of $S$, such that $(IS:x^{e-1})\subset \maxim T$. Hence $a\in\maxim T$. Since $T$ maps to $R^+$, we get $a\in\maxim  R^+$, and hence $a\in  x{\overline{R^+}}$, as we wanted to show.
\qed\end{pf}

The condition that $\maxim$ is cyclic modulo a regular sequence is in this case equivalent with $R$ being \CM\ with regularity defect at most one (recall that the \emph{regularity defect} of $R$ is by definition the difference between its embedding dimension and its Krull dimension). If $R$ is regular, then \eqref{eq:socin} is true for any $\maxim$-primary ideal $I$ of $R$ (use the fact that $R\to  R^+$ is flat).

\subsubsection*{Proof of Theorem~\ref{T:sing}}
Let $(R,\maxim)$ be as in the statement of Theorem~\ref{T:sing}. In particular, $R$ is \CM\ by Theorem~\ref{T:pure}. Let $M$ be a finitely generated $R$-module such that $M/\id M$ has finite length. Let $I$ be the annihilator of $M$. By Nakayama's Lemma, $M/\id M$  having finite length implies that $I+\id$ is $\maxim$-primary. Therefore, we can find a regular sequence $\rij xd$ with $x_1,\dots,x_h\in I$ and $x_{h+1},\dots, x_d\in\id$. We will show by downward induction on $i$ that $R^+/\rij xiR^+$ is flat over $R/\rij xiR$, for all $i\geq h$. In particular, for $i=h$, we get that 
	\begin{equation*}
	\tor {R/\rij xhR}n {R^+/\rij xhR^+}N=0,
	\end{equation*} 
for all $n\geq 1$ and all $R/\rij xhR$-modules $N$.  Since $\rij xh$ is both $R$-regular and $R^+$-regular, we get the required vanishing, by taking $N$ to be $M$. To prove the claim, the case $i=d$ has already been established in the course of proving Theorem~\ref{T:tor}. So assume $h\leq i<d$. In general,  let $A$ be a ring, $x$ an $A$-regular element, $K$ an $A/xA$-module and $L$ an $A$-module. The standard spectral sequence
	\begin{equation*}
	 \tor {A/xA}pK{\tor Aq  {A/xA}L} \implies \tor A{p+q}KL
	\end{equation*}
 degenerates into an exact sequence
	\begin{multline*}
	\tor {A/xA}{n-1}K{\ann  Lx} \to \tor AnKL \to \tor {A/xA}nK{L /xL} \to \dots\\
	K\tensor \ann  Lx \to \tor A1KL\to \tor {A/xA}1K{L/xL} \to 0.
	\end{multline*}
Put $A:=R/\rij xiR$ and   $B:=R^+/\rij xiR^+$, so that by our induction hypothesis $A/x_{i+1}A\to B/x_{i+1}B$ is flat and we need to show that $B$ is flat over $A$. Applying the above spectral sequence with $x:=x_{i+1}$ to the $A/xA$-module $K:=B/xB$ and to an arbitrary $A$-module $L$, we get that $\tor An{B/xB}L=0$, for all $n\geq 2$.   Since $x$ is $B$-regular, the short exact sequence
	\begin{equation*}
	\exactseq B x B {} {B/xB}
	\end{equation*}
gives rise to a long exact sequence
	\begin{equation*}
	\tor A{n+1}{B/xB}L\to \tor AnBL \map x \tor AnBL,
	\end{equation*}
for all $n\geq 1$. Therefore, multiplication with $x$ on $\tor AnBL$ is injective, for all $n\geq 1$. In particular, we have for each $n$ an embedding 
	\begin{equation}\label{eq:emb}
	\tor AnBL\subset (\tor AnBL)_x=\tor {A_x}n{B_x}{L_x}.
	\end{equation} 
Since $x\in\id$, we have that $R_x$ is regular. Therefore, $R_x\to (R_x)^+$ is flat. An easy calculation shows that $(R_x)^+=(R^+)_x$ (see \citep[Lemma 6.5]{HHbigCM}). In particular, $\tor {R_x}n{(R^+)_x}{L_x}=0$, for every $A$-module $L$. Since $\rij xi$ is $R_x$-regular and $(R^+)_x$-regular, we get that $\tor {A_x}n{B_x}{L_x}=0$. Therefore, $\tor AnBL=0$ by \eqref{eq:emb}, for every $A$-module $L$ and every $n\geq 1$, showing that $B$ is flat over $A$, as required. 
\qed

Theorem~\ref{T:sing} implies that  for  $R$ of dimension   three, if $\tor R1{R^+}k$ vanishes, then so does $\tor Rn{R^+}{R/\pr}$ for every $n\geq 1$ and every prime ideal $\pr$ of $R$ not in the singular locus of $R$, since $R$ is normal by Theorem~\ref{T:pure} and hence the ideal defining the singular locus of $R$ has height at least two. On the other hand, we have the following non-vanishing result.

\begin{cor}\label{C:comp}
Let $(R,\maxim)$ be an excellent  local domain   of positive \ch. If $\pr$ is a prime ideal defining an irreducible component of the singular locus of $R$, then $\tor R1{R^+}{R/\pr}$ is non-zero.
\end{cor}
\begin{pf}
Assume $\tor R1{R^+}{R/\pr}$ vanishes. Hence  so does $\tor{R_\pr}1{(R^+)_\pr}{k(\pr)}$, where $k(\pr)$ is the residue field of $\pr$. Since $(R^+)_\pr$ is equal to $(R_\pr)^+$ by \citep[Lemma 6.5]{HHbigCM}, and since $R_\pr$ has an isolated singularity, it follows from Theorem~\ref{T:main} that $R_\pr$ is regular, contradicting the choice of $\pr$.
\qed\end{pf}

In view of Lemma~\ref{L:rig} we can generalize this even further: if $R$ is \CM, then each $\tor Rn{R^+}{R/\pr}$ is non-zero, for $n\geq 1$ and for $\pr$ defining an irreducible component of the singular locus of $R$.

\end{document}